\documentclass[11pt]{article}
\usepackage{amsmath, amssymb, amsthm} 
\usepackage{amsfonts} 
\usepackage{geometry} 
\usepackage{hyperref} 
\usepackage{graphicx} 
\usepackage{enumitem} 

\newtheorem{theorem}{Theorem}[section]
\newtheorem{lemma}[theorem]{Lemma}
\newtheorem{proposition}[theorem]{Proposition}
\newtheorem{corollary}[theorem]{Corollary}
\theoremstyle{definition}
\newtheorem{definition}[theorem]{Definition}
\theoremstyle{remark}

\newtheorem{example}[theorem]{Example}

\newcommand \s{^{*}}
\newcommand \+{^{\dag}}
\newcommand \p{^{\perp}}

\begin{document}

	\title{Convergence of Complementable Operators}
\author{
$^{^{^a}}$Sachin Manjunath Naik and $^{^{^a}}$P. Sam Johnson\\ 	 
	\small $^{^{^a}}$	Department of Mathematical and Computational Sciences,\\
	National Institute of Technology Karnataka, Surathkal,\\
	Mangalore, Karnataka, India, 575025
}
	\date{} 
	
	\maketitle
	
	\begin{abstract}
	Complementable operators extend classical matrix decompositions, such as the Schur complement, to the setting of infinite-dimensional Hilbert spaces, thereby broadening their applicability in various mathematical and physical contexts. This paper focuses on the convergence properties of complementable operators, investigating when the limit of sequence of complementable operators remains complementable. We also explore the convergence of sequences and series of powers of complementable operators, providing new insights into their convergence behavior. Additionally, we examine the conditions under which the set of complementable operators is the subset of set of boundary points of the set of non-complementable operators with respect to the strong operator topology. The paper further explores the topological structure of the subset of complementable operators, offering a characterization of its closed subsets.
	\end{abstract}

\section{Introduction}\label{sec1}
The study of complementable operators has its origins in the analysis of the Schur complement. Complementable operators extend classical matrix decompositions, such as the Schur complement, into the setting of infinite-dimensional Hilbert spaces, thus broadening their scope and applicability. 
The classical Schur complement originated in the study of partitioned matrices. Let $T$ be a block matrix of size $(n+m) \times (n+m)$, partitioned as
$$
T = \begin{pmatrix} A & B \\ C & D \end{pmatrix},
$$
where $A$ is an $n \times n$ matrix, $B$ is an $n \times m$ matrix, $C$ is an $m \times n$ matrix, and $D$ is an $m \times m$ matrix. If the matrix $D$ is nonsingular, the classical Schur complement (denoted by $T_{/D}$) of $D$ in $T$ is given by
$$
T_{/D} = A - BD^{-1}C.
$$
This definition plays a fundamental role in matrix analysis and has applications in various mathematical and physical domains. However, the classical Schur complement refers only to a matrix of lower order derived from a partitioned matrix. In contrast, the generalization of this notion to operators in Hilbert spaces leads to the concept of complementable operators, which are often more complex and useful for infinite-dimensional spaces. Ando \cite{Ando} introduced the concept of complementable operators in $\mathbb{C}^n$ with respect to a subspace $M$ of $\mathbb{C}^n$, extending the classical Schur complement to a more generalized framework. An $ n \times n $ complex matrix $T$ is said to be $ M $-complementable if there exist matrices $ M_r $ and $ M_\ell $ such that $ P M_r = M_r $, $ M_\ell P = M_\ell $, $ P T M_r = P T $ and $ M_\ell T P = T P $, where $ P $ is the orthogonal projection onto $ M $. According to this definition, the operators $ M_\ell T $ and $ T M_r $ are uniquely determined that are independent of the specific choices of $ M_\ell $ and $ M_r $. Further, Ando defined the Schur complement of $T$ with respect to $M$ as  $$ T_{/M} = T - M_\ell T = T - T M_r .$$ It is observed that if $ T $ has the block form $\begin{pmatrix} A & B\\ C & D \end{pmatrix} $ with respect to the subspace $ M $ and if $\mathcal{R}(D)$ is closed, then $$ T_{/M}=\begin{pmatrix} A - B D^{-1} C & 0 \\ 0 & 0 \end{pmatrix}. $$

Ando's framework laid the foundation for various extensions and generalizations.
Mitra \cite{Mitra} extended the idea of complementable operators to the context of operators between Hilbert spaces, introducing the notion of bilateral shorted operators and analyzing their properties.  A comprehensive survey of the Schur complement can be found in \cite{Arias, butler2, Carlson1, carlson3, krien2, Cottle, shorting, Goller, krein3,  Ouellette, Zhang}.  
Using block matrix representation of operators, Antezana et al. \cite{Antezana} introduced the notions of complementability and weak complementability for bounded operators $T: \mathcal{H} \to \mathcal{K}$ with respect to closed subspaces $M \subseteq \mathcal{H}$ and $N \subseteq \mathcal{K}$, based on range inclusion conditions and the modulus of operators. This led to the introduction of the Schur complement (bilateral shorted operator) and its application in the study of parallel sums of bounded linear operators. This Schur complement operator retains key properties from the matrix case.
Arias et al. \cite{Arias} further investigated the relationships between these concepts and the existence of specific bounded projections, providing several characterizations of complementable operators. This definition of the Schur complement introduced by Antezana \cite{Antezana} and Arias et al. \cite{Arias} relies on the notion of weak complementability of operators. Later, this was simplified for complementable operators by Naik and Johnson \cite{ourpaper} and provided a characterization of complementable operators in terms of the images of the unit ball under the components of the block matrix of $T$. These results were used to explore $(M,N)$-complementability and the hypo-$EP$ property of operators $T$. More details of the complementable operators can be found in \cite{ Antezana, Arias, corach2002, Massey}.

In Section $2$, we introduce the basic definitions of complementable operators and provide key characterizations of such operators. Building upon these characterizations, we define $(M,N,\lambda)$-complementability. Additionally, we provide the definition of the Schur complement for a complementable operator. This section also introduces several notations that are utilized throughout the article for ease of reference.

Section 3 is dedicated to the study of the convergence properties of complementable operators. Prior to addressing the convergence of complementable operators specifically, we first present general definitions related to the convergence of operators and explore how these concepts relate to the block matrix representation of operators. We then discuss various convergence properties of complementable operators, including a result showing that a certain class of operators can be expressed as the strong limit of operators that are not $(M,N)$-complementable.

In Section 4, we present some topological results concerning the set of $(M,N,\lambda)$-complementable operators. Specifically, we provide a characterization of the closed subsets of the set of $(M,N)$-complementable operators and examine their topological structure.

\section{Preliminaries}\label{sec2}
Let $\mathcal H$ and $\mathcal K $ be real or complex Hilbert spaces.  We denote  $\mathcal L(\mathcal H, \mathcal K),$  the space of all linear operators from $\mathcal H$ to $\mathcal K.$ We denote  $\mathcal B(\mathcal H, \mathcal K),$  the space of all linear bounded operators from $\mathcal H$ to $\mathcal K$ and we abbreviate $\mathcal B(\mathcal H)=\mathcal B(\mathcal H, \mathcal H)$. For $T \in \mathcal B(\mathcal H, \mathcal K),$ we denote by $T^*$, $\mathcal N(T)$ and $\mathcal R(T),$ respectively, the adjoint, the null-space and the range-space of $T$. For a subspace $M \subseteq \mathcal{H},$  we denote the set $\{x \in  M : \|x\| \leq 1 \}$ as $\mathcal{B}_M$. An operator $T\in \mathcal B(\mathcal H)$ is said to be a projection if $T^2=T.$ A projection, $T$ is said to be orthogonal if $T=T^*.$  An operator $T\in \mathcal B(\mathcal H)$ is said to be self-adjoint if $T\s=T.$ A a self-adjoint operator is called positive (denoted $T\geq 0$) if  $\langle Tx, x \rangle \geq 0$ for all $x\in \mathcal H$. For $T\in \mathcal B(\mathcal H),$ the square root of $T$ is an operator $S \in \mathcal{B}(\mathcal{H})$ such that $S^2 = T$. If $T$ is a positive operator, the square root $S$ is unique and denoted by $T^{1/2}$. The operator $|T|:=(T^*T)^{1/2}$ is called the modulus of $T$. For $T\in \mathcal{B}(\mathcal{H},\mathcal{K})$, the reduced minimum modulus of $T$ is defined as $$\gamma (T)=\inf \{\|Tx\| : x\in \mathcal{H} \cap \mathcal{N}(T){^{\perp}}~ \text{and}~ \|x\|=1 \}.$$ 
Moreover, for a non-zero bounded operator, $\mathcal{R}(T) $ is closed if and only if $\gamma (T)>0$ \cite{minmod}. 
The following theorem, called Douglas' lemma (also known as Douglas' factorization theorem), provides conditions under which a bounded linear operator can be factored.

\begin{theorem}[\cite{Douglas}]\label{dgls}
	Let $ A, B \in \mathcal{B}(\mathcal{H}) $. Then the following are equivalent:
	\begin{enumerate}
		\item $ \mathcal{R}(A) \subseteq \mathcal{R}(B) $;
		\item $ AA^* \leq \lambda BB^* $, for some $ \lambda > 0 $;
		\item There exists a bounded operator $ C \in \mathcal{B}(\mathcal{H}) $ such that $ A = BC $.
	\end{enumerate}
	Moreover, if these equivalent conditions hold, then there is a unique operator $ C \in \mathcal{B}(\mathcal{H}) $ such that
	\begin{enumerate}
		\item[(i)] $ \|C\| = \inf \{ \lambda > 0 : AA^* \leq \lambda BB^* \} $;
		\item[(ii)] $ \mathcal{N}(A) = \mathcal{N}(B) $;
		\item[(iii)] $ \mathcal{R}(C) \subseteq \mathcal{N}(B)^\perp $.
	\end{enumerate}
\end{theorem}
\noindent The unique solution $ C $ is referred as the Douglas reduced solution \cite{redsoln}.

Let $ T: \mathcal{H} \to \mathcal{K} $ be a linear operator and let $ M $ and $ N $ be closed subspaces of $\mathcal{H}$ and $\mathcal{K}$ respectively. We can express $\mathcal{H}$ as the orthogonal direct sum of $ M $ and $ M^\perp $ and $\mathcal{K}$ as the orthogonal direct sum of $ N $ and $ N^\perp $. With respect to the above decompositions, the operator $ T $ can be written in a block matrix form as
\begin{equation}\label{*}
	T = \begin{pmatrix}
		A & B \\
		C & D
	\end{pmatrix}.
\end{equation}

\begin{definition}[\cite{Antezana}]
	Let $ P_r \in \mathcal{B}(\mathcal{H}) $ and $ P_\ell \in \mathcal{B}(\mathcal{K}) $ be projections. An operator $ T \in \mathcal{B}(\mathcal{H}, \mathcal{K}) $ is said to be $(P_r, P_\ell)$-complementable if there exist operators $ M_r \in \mathcal{B}(\mathcal{H}) $ and $ M_\ell \in \mathcal{B}(\mathcal{K}) $ such that
	\begin{enumerate}
		\item $ (I-P_r)M_r=M_r $ and $ (I-P_\ell)TM_r=(I-P_\ell) T $;
		\item $ M_\ell (I-P_\ell)=M_\ell $ and $ M_\ell T(I-P_r)=T(I-P_r) $.
	\end{enumerate}
\end{definition}

The above definition is independent of the specific projections considered and it depends solely on the ranges of these projections. The following theorem demonstrates this dependency on the ranges of the projections.

\begin{proposition}[\cite{Antezana}] \label{propdef}
	Let $ P_r \in \mathcal{B}(\mathcal{H}) $ and $ P_\ell \in \mathcal{B}(\mathcal{K}) $ be projections whose ranges are $ M $ and $ N $ respectively. Let $ T \in \mathcal{B}(\mathcal{H}, \mathcal{K}) $ be as given in (\ref{*}), the following statements are equivalent:  
	\begin{enumerate}
		\item $ T $ is $(P_r, P_\ell)$-complementable;
		\item $ \mathcal{R}(C) \subseteq \mathcal{R}(D) $ and $ \mathcal{R}(B^*) \subseteq \mathcal{R}(D^*) $;
		\item There exist two projections $ \tilde{P} \in \mathcal{B}(\mathcal{H}) $ and $ \tilde{Q} \in \mathcal{B}(\mathcal{K}) $ such that: $ \mathcal{R}(\tilde{P}^*)=M $, $ \mathcal{R}(\tilde{Q})=N $,  $ \mathcal{R}(T\tilde{P})\subseteq N $ and $ \mathcal{R}((\tilde{Q}T)^*)\subseteq M $.
	\end{enumerate}
\end{proposition}

\begin{definition}[\cite{Antezana}]
	An operator \quad $ T \in \mathcal{B}(\mathcal{H},  \mathcal{K}) $ \quad  is said to be \quad $(M, N)$-complementable if it is $(P_r, P_\ell)$-complementable for some projections $ P_r $ and $ P_\ell $ with $ \mathcal{R}(P_r)=M $ and $ \mathcal{R}(P_\ell)=N $.
\end{definition}
\begin{theorem}(\cite{ourpaper})\label{char2}
	Let  $T \in \mathcal{B}(\mathcal{H},\mathcal{K})$ be as given in (\ref{*}). Then $T$ is $(M,N)$-complementable if and only if there exists $\lambda > 0$ such that  \begin{equation}\label{lambda}
		C(\mathcal{B}_{M}) \subseteq \lambda D(\mathcal{B}_{ M\p})~ \text{and}~ B\s(\mathcal{B}_{N}) \subseteq \lambda D\s(\mathcal{B}_{ N\p}).
	\end{equation}
	
	Moreover, if $\mathcal{R}(D)$ is a non-zero closed subspace of $N\p$, then $T$ is $(M,N)$-\\
	complementable if and only if $$C(\mathcal{B}_{M}) \subseteq \frac{\|C\|}{\gamma (D)}D(\mathcal{B}_{ M\p}) \quad \text{and} \quad B\s(\mathcal{B}_{N}) \subseteq \frac{\|B\|}{\gamma (D\s)}D\s(\mathcal{B}_{N\p})$$
\end{theorem}
\noindent If (\ref{lambda}) holds, then $T$ is said to be $(M,N, \lambda)$-complementable. 

In finite-dimensional Hilbert spaces, every positive operator $ T $ is $ M $-complementable. Indeed, if $T$ is the positive operator of the form (\ref{*}), then the inclusions $ \mathcal{R}(C) \subseteq \mathcal{R}(D) $  and $ \mathcal{R}(B\s) \subseteq \mathcal{R}(D\s) $ always hold. The following example shows that every positive operator is not $M$-complementable in infinite dimensional space.
\begin{example}
	Let $D: \ell_2 \to \ell_2$ be defined by $$D(x_1,x_2,x_3,\ldots)=\Big(x_1,\frac{x_2}{4},\frac{x_3}{9},\ldots\Big).$$
	Let $$M=\ell_2 \oplus 0 ~\text{and} ~ M_0= \text{span}\Big(\Big\{\frac{1}{n}\Big\}\Big).$$
	Now consider $T:\ell_2 \oplus \ell_2 \to \ell_2 \oplus \ell_2 $ defined by $$T=\begin{pmatrix}
		D & P_{M_0} \\
		P_{M_0} & D
	\end{pmatrix}.$$
	
\end{example}
Here, we have  $\mathcal{R}(P_{M_0}) \nsubseteq \mathcal{R}(D).$ So, $T$ is not $M$-complementable. \\ Now, we have 
\begin{align*}
	\langle T(x+y), x+y \rangle &= \langle Dx,x \rangle +\langle P_{M{_0}} y,x \rangle + \langle y ,P_{M{_0}} x \rangle + \langle y,Dy \rangle\\
	&= \langle Dx,x \rangle +\langle  y,P_{M{_0}} x \rangle + \langle y ,P_{M{_0}} x \rangle + \langle y,Dy \rangle.
\end{align*} Now $y= \alpha \{\frac{1}{n}\} + z$, for some  $ z \in ({M{_0}}) \p$, $\alpha >0$ and $P_{M{_0}} x = \beta \{\frac{1}{n}\}$ for some $ \beta >0.$
Thus, $$\langle T(x+y), x+y \rangle= \langle Dx,x \rangle + \langle y,Dy \rangle + 2 \alpha \beta \Big\|\Big\{\frac{1}{n}\Big\}\Big \|^2 \geq 0.$$
Hence $T$ is a positive operator.

Motivated by this observation, Antezana et al. \cite{Antezana} defined a weaker notion of complementable operator called weak complementable operator as follows:

\begin{definition}[\cite{Antezana}] 
	An operator $ T \in \mathcal{B}(\mathcal{H}, \mathcal{K}) $ of the from (\ref{*}) is said to be $(M,N)$-weakly complementable if
	$$ \mathcal{R}(C) \subseteq \mathcal{R}(|D^*|^{\frac{1}{2}}) \ \text{and} \ \mathcal{R}(B^*) \subseteq \mathcal{R}(|D|^{\frac{1}{2}}).$$
\end{definition}
It is easy to see that if $ T $ is $(M,N)$-complementable, then it is $(M,N)$-weakly complementable. Both notions of complementability coincide if $ \mathcal{R}(D) $ is closed.

\begin{definition}[\cite{Antezana}]
	Let $ T $ be $(M,N)$-weakly complementable operator of the form (\ref{*}). Let $ F $ and $ E $ be Douglas reduced solutions of the equations $ C=|D^*|^{\frac{1}{2}}UX $ and $ B^* = |D|^{\frac{1}{2}}X $ respectively, where $ D=U|D| $ is the polar decomposition of $ D $. Then the Schur complement (bilateral shorted operator) of $ T $ with respect to the closed subspaces $ M $ and $ N $ is
	$$  T_{/(M,N)}=\begin{pmatrix}
		A-E^*F & 0 \\
		0 & 0
	\end{pmatrix}. $$
	If $ \mathcal{H}=\mathcal{K} $ and $ M=N $, then we denote $ T_{/(M,M)} $ by $ T_{/M} $, which is said to be the Schur complement of $ T $ with respect to $ M $.
\end{definition}
\noindent Moreover, it is observed in \cite{ourpaper} that if $T$ is $(M,N)$-complementable, then the Schur complement of $T$ with respect to the closed subspaces $M$ and $N$ can be written as $$T_{/(M,N)}=\begin{pmatrix}
	A-BZ & 0 \\
	0 & 0
\end{pmatrix}=\begin{pmatrix}
	A-YC & 0 \\
	0 & 0
\end{pmatrix},$$ where $Z$ is the Douglas reduced solution of $\mathcal{R}(C) \subseteq \mathcal{R}(D)$ and $Y\s$ is the Douglas reduced solution of $\mathcal{R}(B\s) \subseteq \mathcal{R}(D\s).$\\
For a given $\lambda > 0$, we denote 
$$
\psi(M,N,\lambda) = \{ T \in \mathcal{B}(\mathcal{H}, \mathcal{K}) : T \quad \text{is} \quad (M,N,\lambda)\text{-complementable} \}.
$$
Now for $\lambda_1 < \lambda_2$, we have $\psi(M,N,\lambda_1) \subsetneq \psi(M,N,\lambda_2)$.\\
We denote the set of all $(M,N)$-complementable operators by $\psi(M,N)$. Thus, by our notation, it follows that 
$$
\psi(M,N) = \bigcup_{\lambda \in \mathbb{N}} \psi(M,N,\lambda).
$$
We use $\mathcal{C}_D(M,N)$ to denote the set of operators $T \in \mathcal{L}(\mathcal{H}, \mathcal{K})$ with $\mathcal{R}(D)$ is closed and $\mathcal{B}\mathcal{C}_D(M,N)$ to denote the set of operators $T \in \mathcal{B}(\mathcal{H}, \mathcal{K})$ with $\mathcal{R}(D)$ is closed when $T$ is written in the form (\ref{*}).

\section{Convergence of Complementable Operators}

The behavior of bounded operators as they converge towards a limit is a crucial area of investigation in Hilbert spaces. The notion of convergence enables a rigorous analysis of operator behavior, facilitating the study of their stability, spectral properties, and various applications in both theoretical and applied mathematics. Understanding how operators behave under limiting processes helps in exploring questions related to their continuity, the influence of small perturbations, and how they approximate each other in different contexts.

In the specific case of complementable operators, convergence is particularly important when analyzing how operators behave in relation to their complements or when approximating operators in the context of subspaces. Before diving into the convergence of complementable operators, it is beneficial to first define and explore general concepts of convergence in operator theory. This will provide a solid framework for understanding the nuances of convergence in the block matrix representation of operators.

\begin{definition}[\cite{kreyzig}]
	Let $X$ and $Y$ be normed spaces. Let $\{T_n\}$ be a  sequence operators in $ \mathcal{B}(X,Y).$
	\begin{enumerate}
		\item A sequence $\{T_n\}$ of operators in $ \mathcal{B}(X,Y)$ is said to be  uniformly convergent if $\{T_n\}$ converges in the norm on $\mathcal{B}(X,Y).$
		\item A sequence $\{T_n\}$ of operators in $ \mathcal{B}(X,Y)$ is said to be is said to be strongly convergent if $\{T_n x\}$ converges in $Y$
		for every $x \in X.$
	\end{enumerate}
\end{definition}
\noindent If $\{T_n\}$ is uniformly convergent to $T$, we denote it by $T_n \to T$ and 

\noindent If $\{T_n\}$ is strongly convergent to $T$, we denote it by $ T_n \xrightarrow{\text{s}} T$.

The following theorem provides a criterion for the convergence of a sequence of operators $\{T_n\}$ in $\mathcal{B}(\mathcal{H},\mathcal{K})$, where

\begin{equation}\label{Tn}
	T_n = \begin{pmatrix}
		A_n & B_n \\
		C_n & D_n
	\end{pmatrix}
\end{equation}
to an operator
\begin{equation}\label{T}
	T = \begin{pmatrix}
		A & B \\
		C & D
	\end{pmatrix}
\end{equation}
in $\mathcal{B}(\mathcal{H},\mathcal{K})$ with respect to the convergence of the sequence of the operators  $\{A_n\}$, $\{B_n\}$, $\{C_n\}$ and $\{D_n\}.$
\begin{theorem}\label{thmblock}
	Let $\{T_n\} \subseteq \mathcal{B}(\mathcal{H}, \mathcal{K})$ and $T \in \mathcal{B}(\mathcal{H}, \mathcal{K})$ be as defined in (\ref{Tn}) and (\ref{T}) respectively. The sequence $\{T_n\}$ converges uniformly to $T$ if and only if the components $\{A_n\}$, $\{B_n\}$, $\{C_n\}$ and $\{D_n\}$  converge uniformly to $A$, $B$, $C$ and $D$ respectively.
	
\end{theorem}

\begin{proof}
	Let $P_M$, $P_{M\p}$, $P_N$ and $P_{N\p}$ be orthogonal projections onto the subspaces $M$, $M\p$, $N$ and $N\p$ respectively. Now we have 
	\begin{align}
		P_N T_n P_M &= \begin{pmatrix}
			A_n & 0 \\
			0 & 0
		\end{pmatrix} \quad \text{and} \quad P_N T P_M = \begin{pmatrix}
			A & 0 \\
			0 & 0
		\end{pmatrix}, \\
		P_N T_n P_{M\p} &= \begin{pmatrix}
			0 & B_n \\
			0 & 0
		\end{pmatrix} \quad \text{and} \quad P_N T P_{M\p} = \begin{pmatrix}
			0 & B \\
			0 & 0
		\end{pmatrix}, \\
		P_{N\p} T_n P_M &= \begin{pmatrix}
			0 & 0 \\
			C_n & 0
		\end{pmatrix} \quad \text{and} \quad P_{N\p} T P_M = \begin{pmatrix}
			0 & 0 \\
			C & 0
		\end{pmatrix}, \\
		P_{N\p} T_n P_{M\p} &= \begin{pmatrix}
			0 & 0 \\
			0 & D_n
		\end{pmatrix} \quad \text{and} \quad P_{N\p} T P_{M\p} = \begin{pmatrix}
			0 & 0 \\
			0 & D
		\end{pmatrix}.
	\end{align}
	Now, it is enough to show that the sequence $\{T_n\}$ converges uniformly to $T$ if and only if $\{P_N T_n P_M\}$, $\{P_N T_n P_{M\p}\}$, $\{P_{N\p} T_n P_M\}$ and $\{P_{N\p} T_n P_{M\p}\}$ converges uniformly  to $P_N T P_M$, $P_N T P_{M\p}$, $P_{N\p} T P_M$ and $P_{N\p} T P_{M\p}$ respectively.
	
	\noindent \textbf{($\implies$):} Assume that $\{T_n\}$ converges uniformly to $T$. Hence, for a given $\varepsilon > 0$, there exists $n_0 \in \mathbb{N}$ such that $\|T_n - T\| \leq \varepsilon$ for all $n > n_0.$ \\
	Therefore,
	$$
	\|P_N T_n P_M - P_N T P_M\|  \leq \|P_N\|\|T_n - T\|\|P_M\| \leq \|T_n - T\| \leq \varepsilon.
	$$
	So, $$\{P_N T_n P_M\} \text{ converges uniformly  to } P_N T P_M.$$ 
	Similarly, $\{P_N T_n P_{M\p}\}$, $\{P_{N\p} T_n P_M\}$ and $\{P_{N\p} T_n P_{M\p}\}$ converge uniformly  to $P_N T P_{M\p}$, $P_{N\p} T P_M$ and $P_{N\p} T P_{M\p}$ respectively.
	
	\noindent \textbf{($\impliedby$):} Conversely, assume that $\{P_N T_n P_M\}$, $\{P_N T_n P_{M\p}\}$, $\{P_{N\p} T_n P_M\}$ and $\{P_{N\p} T_n P_{M\p}\}$ converge uniformly to $P_N T P_M$, $P_N T P_{M\p}$, $P_{N\p} T P_M$ and $P_{N\p} T P_{M\p}$ respectively.\\
	Hence, for a given $\varepsilon > 0$, there exists $n_0 \in \mathbb{N}$ such that for all $n > n_0$, 
	\begin{align*}
		&\|P_N T_n P_M - P_N T P_M\| ~~~~~~~< \varepsilon, \\
		&\|P_N T_n P_{M\p} - P_N T P_{M\p}\| ~~~< \varepsilon, \\
		&\|P_{N\p} T_n P_M - P_{N\p} T P_M\| ~~~< \varepsilon, \\
		&\|P_{N\p} T_n P_{M\p} - P_{N\p} T P_{M\p}\| < \varepsilon.
	\end{align*}
	We have 
	\begin{equation*}
		T_n = P_N T_n P_M + P_N T_n P_{M\p} + P_{N\p} T_n P_M + P_{N\p} T_n P_{M\p}
	\end{equation*}
	and 
	\begin{equation*}\label{Tsum}
		T = P_N T P_M + P_N T P_{M\p} + P_{N\p} T P_M + P_{N\p} T P_{M\p}.
	\end{equation*}
	So, 
	\begin{align*}
		T_n - T &= (P_N T_n P_M + P_N T_n P_{M\p} + P_{N\p} T_n P_M + P_{N\p} T_n P_{M\p}) \\
		&~~~~- (P_N T P_M + P_N T P_{M\p} + P_{N\p} T P_M + P_{N\p} T P_{M\p}) \\
		&= P_N (T_n - T) P_M + P_N (T_n - T) P_{M\p} \\
		&~~~+ P_{N\p} (T_n - T) P_M + P_{N\p} (T_n - T) P_{M\p}.
	\end{align*}
	Thus, 
	\begin{align*}
		\|T_n - T\| &\leq \|P_N\| \|T_n - T\| \|P_M\| + \|P_N\| \|T_n - T\| \|P_{M\p}\| \\
		&\quad + \|P_{N\p}\| \|T_n - T\| \|P_M\| + \|P_{N\p}\| \|T_n - T\| \|P_{M\p}\| \\
		&\leq \|T_n - T\| + \|T_n - T\| + \|T_n - T\| + \|T_n - T\|
		\leq 4\varepsilon.
	\end{align*}
	Hence $\{T_n\}$ converges uniformly to $T$.
\end{proof}

\begin{theorem}
	Let $T \in \mathcal{B}(\mathcal{H}, \mathcal{K})$ be as defined in (\ref{T}). Then $$\max \{\|A\|, \|B\|, \|C\|, \|D\| \} \leq \|T\| \leq \|A\|+ \|B\| + \|C\| + \|D\|.$$
\end{theorem}

\begin{proof}
	Let $x \in M$. Then we have 
	$$
	Tx = \begin{pmatrix}
		A & B\\
		C & D
	\end{pmatrix}\begin{pmatrix}
		x \\ 0
	\end{pmatrix} = \begin{pmatrix}
		Ax \\ Cx
	\end{pmatrix}.
	$$
	So, we have 
	$$
	\|Tx\|^2 = \|Ax\|^2 + \|Cx\|^2 \geq \|Ax\|^2.
	$$
	Thus, $\|Ax\| \leq \|Tx\|$ for all $x \in M$. Therefore, we have $\|A\| \leq \|T\|$. Similarly, we can prove $\|B\| \leq \|T\|$, $\|C\| \leq \|T\|$ and $\|D\| \leq \|T\|$. \\
	Hence,
	\begin{equation}\label{lineq}
		\max \{\|A\|, \|B\|, \|C\|, \|D\|\} \leq \|T\|.
	\end{equation}
	Also,  
	\begin{align*}
		\|T\| &= \|P_N T P_M + P_N T P_{M\p} + P_{N\p} T P_M + P_{N\p} T P_{M\p}\| \\
		&\leq \|P_N T P_M\| + \|P_N T P_{M\p}\| + \|P_{N\p} T P_M\| + \|P_{N\p} T P_{M\p}\|.
	\end{align*}
	Hence,
	\begin{equation}\label{rineq}
		\|T\| \leq \|A\| + \|B\| + \|C\| + \|D\|.
	\end{equation}
	Thus,
	$$\max \{\|A\|, \|B\|, \|C\|, \|D\| \} \leq \|T\| \leq \|A\|+ \|B\| + \|C\| + \|D\|.$$
\end{proof}

Now we discuss the convergence of  $(M,N)$-complementable operators. This property of  $(M,N)$-complementable operators are vital in operator theory, particularly in applications where complementability ensures the stability of certain operator properties under approximation. One of the fundamental questions is whether the limit sequence of complementable operator is again a complementable operator or not? One can ask this question with respect to various operator topologies, such as the norm or strong operator topology. In general, it is not closed under the norm topology. The following example demonstrates that the limit of a sequence of $(M,N)$-complementable operators may fail to retain this property.

\begin{example}  
	Let $\mathcal{H} = \ell_2 \oplus \ell_2$ and $M = \ell_2 \oplus \{0\}$. For each $n \in \mathbb{N},$ define the operators $T_n \in \mathcal{B}(\mathcal{H})$ with respect to the decomposition $\mathcal{H} = \ell_2 \oplus \ell_2$ as:  
	$$
	T_n = \begin{pmatrix}
		A_n & B_n \\
		C_n & D_n
	\end{pmatrix},
	$$  
	where $A_n, B_n, C_n,$ and $D_n$ are given by 
	$$
	A_n = I, \quad B_n = I, \quad C_n = \left(1 + \frac{1}{n}\right)I, \quad D_n = \frac{1}{n}I, \text{ for all } n \in \mathbb{N}.
	$$  
\end{example} 
Here, both $C_n$ and $D_n$ are invertible. Since $\mathcal{R}(C_n) \subseteq \mathcal{R}(D_n)$, it follows that each $T_n$ is $M$-complementable.  \\
Now we have that  $$T_n \to T = \begin{pmatrix}
	I & I \\
	I & 0
\end{pmatrix} $$
and  $\mathcal{R}(I) \nsubseteq \mathcal{R}(0)$, the operator $T$ is not $M$-complementable.  

While the previous example shows that the set of $(M,N)$-complementable operators is not necessarily closed, it is still possible to identify conditions under which a sequence of such operators converges to a complementable operator. By imposing specific conditions or certain structural properties on the operators, we can ensure that complementability is preserved under limits. The following result establishes sufficient conditions under which a sequence of $(M,N)$-complementable operators converges to a complementable operator, providing insights into the stability of complementability under approximation.

\begin{theorem}\label{conv}
	Let $\{T_n\} \subseteq \psi(M,N,\lambda_n)$ and $T \in \mathcal{B}\mathcal{C}_D(M,N)$ be as defined in (\ref{Tn}) and (\ref{T}) with $T_n \to T$. If $\lambda_n \|D - D_n\| \to 0$, then $T \in \psi (M,N). $
\end{theorem}

\begin{proof}
	Since $T_n \in \psi(M,N,\lambda_n)$, we have
	$$C_n(\mathcal{B}_{M}) \subseteq \lambda_n D_n(\mathcal{B}_{M\p})~\text{and}~ B_n\s (\mathcal{B}_{N}) \subseteq \lambda_n D_n\s(\mathcal{B}_{N\p}).$$
	Let $\beta_n = \sup \{\lambda_i : 1 \leq i \leq n\}$. Now we can write $C_n(\mathcal{B}_{M}) \subseteq \beta_n D_n(\mathcal{B}_{M\p})$ and $B_n\s(\mathcal{B}_{N}) \subseteq \beta_n D_n\s(\mathcal{B}_{N\p})$. Additionally, $\{ \beta_n \}$ constitutes a monotonically increasing sequence of positive real numbers. Let $x \in \mathcal{B}_{M}$. Since $C_n \to C$, we observe that $C_n x \to Cx$. Consequently, we have $$C_n x = \beta_n D_n y_n \quad \text{for some} \quad y_n \in \mathcal{B}_{M\p} \quad \text{and} \quad \beta_n D_n y_n \to Cx.$$
	We have 
	$$\begin{aligned}
		&\| \beta_n Dy_n - \beta_m Dy_m\| \\
		&= \| \beta_n Dy_n - \beta_n D_n y_n + \beta_n D_n y_n - \beta_m D_m y_m + \beta_m D_m y_m - \beta_m Dy_m\| \\
		&\leq \| (D - D_n)(\beta_n y_n)\| + \| \beta_n D_n y_n - \beta_m D_m y_m\| + \| (D_m - D)(\beta_m y_m)\| \\
		&\leq \|\beta_n\| \|(D - D_n)\| + \| \beta_n D_n y_n - \beta_m D_m y_m\| + \|\beta_m\| \|(D_m - D)\|.
	\end{aligned}
	$$
	Since $\lambda_n \|D - D_n\| \to 0$, it follows that $\beta_n \|D - D_n\| \to 0$ and also $\{\beta_n D_n y_n\}$ is a Cauchy sequence. Thus, for a given $\varepsilon > 0$, there exists $n_0 \in \mathbb{N}$ such that $\|\beta_n\| \|D - D_n\| < \varepsilon$ and $\| \beta_n D_n y_n - \beta_m D_m y_m \| < \varepsilon$ for all $n, m > n_0$. 
	Consequently, we have 
	$$
	\| \beta_n Dy_n - \beta_m Dy_m \| < 3\varepsilon \quad \text{for all} \quad n, m > n_0.
	$$ 
	Therefore, $\{\beta_n Dy_n\}$ converges.\\
	Furthermore, $$\|\beta_n Dy_n - \beta_n D_n y_n\| \leq \beta_n \|D - D_n\| \|y_n\| \leq \beta_n \|D - D_n\| \to 0.$$ This implies that $\{\beta_n Dy_n\}$ converges to the same limit as $\{\beta_n D_n y_n\}$. 
	Thus, $$\beta_n Dy_n = D(\beta_n y_n) \to Cx.$$ This implies $Cx \in \overline{\mathcal{R}(D)}$. Since $\mathcal{R}(D)$ is closed, we conclude that $Cx \in \mathcal{R}(D)$.
	
	Now, let $w \in M$. Since $\frac{w}{\|w\|} \in \mathcal{B}_M,$ we have $C\Big(\frac{w}{\|w\|}\Big)=\frac{Cw}{\|w\|} \in \mathcal{R}(D)$. This implies $Cw \in \mathcal{R}(D)$. Therefore, we have established $\mathcal{R}(C) \subseteq \mathcal{R}(D)$.\\
	By a similar argument, we deduce $\mathcal{R}(B\s) \subseteq \mathcal{R}(D\s)$.
	Hence, $T \in \psi(M,N).$
\end{proof}

\begin{corollary}\label{boundedlambda}
	Let $\{T_n\} \subseteq \psi(M,N,\lambda_n)$ and $T \in \mathcal{B}\mathcal{C}_D(M,N)$ be as defined in (\ref{Tn}) and (\ref{T}) with $T_n \to T$. If $\{\lambda_n\}$ is a bounded sequence of positive real numbers, then $T \in \psi(M,N,\lambda),$ where $\lambda = \sup\{\lambda_n : n \in \mathbb{N}\}$.
\end{corollary}
\begin{proof}
	Given that $\{\lambda_n\}$ is a bounded sequence and $D_n \to D$, we can establish that $\lambda_n \|D_n - D\| \to 0$. Let $\lambda = \sup\{\lambda_n : n \in \mathbb{N}\}$. Thus, we have $$|\lambda_n | \|D_n -D \| \leq |\lambda | \|D_n -D \| \to 0.$$
	Hence, by Theorem \ref{conv}, we have $T$ is $(M,N)$-complementable. Thus $$Cx= \beta D y= \lambda D\Big(\frac{\beta}{\lambda}y\Big), \text{~for some~} y \in {M\p} \cap (\mathcal{N}(D))\p, \beta >0.$$  Hence, we have \begin{equation}\label{bddlambda}
		D(\lambda y_n) \to \lambda D\Big(\frac{\beta}{\lambda}y\Big).
	\end{equation}
	Let $x \in \mathcal{B}_{M}$. Since $C_n \to C$, we have that $C_n x \to Cx$. Since each $T_n$ is $(M,N, \lambda)$-complementable, we have $$C_n x = \lambda D_n y_n \to Cx.$$ By replicating the proof of Theorem \ref{conv}, we can have $$D(\lambda y_n) \to Cx, \quad \text{where} \quad y_n \in \mathcal{B}_{M\p}.$$ 
	Now, for each $n \in \mathbb{N}$, we write $y_n = u_n \oplus v_n,$ for some $u_n \in (\mathcal{N}(D))^\perp$ and $v_n \in \mathcal{N}(D)$ with $\|u_n\| \leq \|y_n\| \leq 1$. Since $D\+ D$ is a projection onto $(\mathcal{N}(D))^\perp$, we have $D\+ D y_n = u_n.$ As $\mathcal{R}(D)$ is closed, $D\+$ is bounded. From this, we can write $D\+ D(\lambda y_n) \to \beta D\+ D y.$ Thus, 
	$$
	\lambda u_n \to \beta y.
	$$
	Since $\|u_n\| \leq 1$ for all $n$, we have 
	$$
	\frac{\beta}{\lambda}\|y\| \leq 1.
	$$
	Therefore, from equation (\ref{bddlambda}), we conclude that 
	$$
	D(\lambda y_n) \to \lambda D(z), \text{ for some } z \in \mathcal{B}_{N^\perp}.
	$$
	In other words, we can write 
	$$
	C(\mathcal{B}_{M}) \subseteq \lambda D(\mathcal{B}_{M^\perp}).
	$$
	
	\noindent Similarly, we can show that $$B\s(\mathcal{B}_{N}) \subseteq \lambda D\s(\mathcal{B}_{N\p}).$$
	Therefore, $T$ is $(M,N, \lambda)$-complementable. 

\end{proof}


\begin{corollary}
	Let $\{T_n\} \subseteq \psi(M,N,\lambda_n)$ and $T \in \mathcal{B}\mathcal{C}_D(M,N)$ be as defined in (\ref{Tn}) and (\ref{T}) with $T_n \to T$. If $ \lambda_n \|D - D_n\| \to 0 $, then $T \in \psi(M,N,\lambda)$ for some $ \lambda $.
\end{corollary}

\begin{proof}
	The proof follows from Theorems \ref{char2} and \ref{conv}.
\end{proof}

\begin{corollary}
	Let $\{T_n\} \subseteq \psi(M,N,\lambda_n)$ and $T_n,T \in \mathcal{B}\mathcal{C}_D(M,N)$ be as defined in (\ref{Tn}) and (\ref{T}) with $T_n \to T $. If $ \frac{\|D - D_n\|}{\gamma(D_n)} \to 0 $, then $T \in \psi(M,N).$
\end{corollary}

\begin{proof}
	By Theorem \ref{char2}, we have $ C(\mathcal{B}_{M}) \subseteq \frac{\|C\|}{\gamma(D_n)} D(\mathcal{B}_{M\p}) $ and $ B\s(\mathcal{B}_{N}) \subseteq \frac{\|B\|}{\gamma(D_n\s)} D\s(\mathcal{B}_{N\p}) $. Since $ \frac{\|D - D_n\|}{\gamma(D_n)} \to 0 $, it follows from Theorem \ref{conv} that $ T $ is $(M,N)$-complementable.
\end{proof}

For a given $T \in \mathcal{B}(\mathcal{H})$, we explore the complementability of powers of sequences and series of $T$. We introduce a sub-collection of $\psi(M, N)$ in which sequences and series of powers of each operator are again complementable operators.
Let us consider $T \in \mathcal{B}(\mathcal{H})$ satisfying 
\begin{equation}\label{A}
	(A + C)(\mathcal{B}_M) \subseteq \lambda (B + D)(\mathcal{B}_{M^\perp}),
\end{equation}
and
\begin{equation}\label{B}
	(A^* + B^*)(\mathcal{B}_M) \subseteq \lambda (C^* + D^*)(\mathcal{B}_{M^\perp}).
\end{equation}
For $x \in \mathcal{B}_M$, we have 
$$
(A + C)x = \lambda (B + D)y,  \text{ for some }  y.
$$
This implies 
$$
Ax + Cx = \lambda By + \lambda Dy.
$$
Rewriting, we get 
$$
Ax - \lambda By = \lambda Dy - Cx.
$$
Note that $Ax - \lambda By \in M$ and $\lambda Dy - Cx \in M^\perp$. Thus, we have 
$$
Ax = \lambda By \quad \text{and} \quad Cx = \lambda Dy.
$$
This gives 
$$
C(\mathcal{B}_M) \subseteq D(\mathcal{B}_{M^\perp}) \quad \text{and} \quad B(\mathcal{B}_{M^\perp}) \subseteq A(\mathcal{B}_M).
$$
Similarly, from  $(A^* + B^*)(\mathcal{B}_M) \subseteq \lambda (C^* + D^*)(\mathcal{B}_{M^\perp})$, we can prove that $$
B^*(\mathcal{B}_M) \subseteq D^*(\mathcal{B}_{M^\perp}) \quad \text{and} \quad C^*(\mathcal{B}_{M^\perp}) \subseteq A^*(\mathcal{B}_M).
$$
Thus, we have 
$$
T \in \psi(M, N, \lambda) \cap \psi(M^\perp, N^\perp, \lambda).
$$

\noindent Now, we introduce the following notations:
\begin{enumerate}
	\item  $T \in \mathcal{B}(\mathcal{H})$ satisfying (\ref{A}) is denoted by $\phi_L(M, N, \lambda)$.
	\item $T \in \mathcal{B}(\mathcal{H})$ satisfying (\ref{B}) is denoted by $\phi_R(M, N, \lambda)$.
	\item $T \in \mathcal{B}(\mathcal{H})$ satisfying both (\ref{A}) and (\ref{B}) is denoted by $\phi(M, N, \lambda)$.
\end{enumerate}
With these notations, we conclude that 
$$
\phi(M, N, \lambda) \subseteq \psi(M, N, \lambda) \cap \psi(M^\perp, N^\perp, \lambda).
$$

\begin{proposition}\label{lem25}
	Let 
	$$
	T_i = \begin{pmatrix}
		A_i & B_i \\
		C_i & D_i
	\end{pmatrix} \in \mathcal{B}(\mathcal{H})
	$$
	for $i=1,2$ and let $\lambda > 0$.
	\begin{enumerate}
		\item If $T_2 \in \phi_L (M,N,\lambda),$  then $T_1T_2 \in \phi_L (M,N,\lambda).$
		\item If $T_1 \in \phi_R (M,N,\lambda),$ then $T_1T_2 \in \phi_R (M,N,\lambda).$
		
		\item If  $T_1 \in \phi_R (M,N,\lambda)$ and $T_2 \in \phi_L (M,N,\lambda),$ then $T_1T_2 \in \phi (M,N,\lambda).$
		\item If  $T_1 \in \phi_R (M,N,\lambda)$ and $T_2 \in \phi_L (M,N,\lambda),$ then $T_1T_2 \in \psi(M,N,\lambda) \cap \psi (M\p,N\p, \lambda).$
	\end{enumerate}
\end{proposition}

\begin{proof}
	Expanding the product $T_1 T_2$, we have
	$$
	T_1 T_2 = \begin{pmatrix}
		A_1 & B_1 \\
		C_1 & D_1
	\end{pmatrix} \begin{pmatrix}
		A_2 & B_2 \\
		C_2 & D_2
	\end{pmatrix} = \begin{pmatrix}
		A_1 A_2 + B_1 C_2 & A_1 B_2 + B_1 D_2 \\
		C_1 A_2 + D_1 C_2 & C_1 B_2 + D_1 D_2
	\end{pmatrix}.
	$$
	
	\noindent (1.) Assume that $T_2 \in \phi_L (M,N,\lambda).$ 
	By assumption, $$(A_2 + C_2)x \in \lambda (B_2 + D_2)(\mathcal{B}_{M^\perp}).$$ 
	Therefore, there exists $y \in \mathcal{B}_{M^\perp}$ such that 
	$$
	(A_2 + C_2)x = \lambda (B_2 + D_2)y.
	$$
	Thus, $$
	A_2 x - B_2 (\lambda y) =  D_2 (\lambda y) - C_2 x.
	$$
	We have that $A_2 x - B_2 (\lambda y) \in M$ and $D_2 (\lambda y) - C_2 x \in M\p.$
	This leads to
	$$
	A_2 x - B_2 (\lambda y) = 0 \quad \text{and} \quad D_2 (\lambda y) - C_2 x = 0.
	$$
	Consequently, we have $$A_2 x = \lambda B_2 y \quad \text{and} \quad C_2 x = \lambda D_2 y.$$
	
	\noindent For any $x \in \mathcal{B}_M$, we have
	\begin{align*}
		((A_1 A_2 + B_1 C_2)&+(C_1 A_2 + D_1 C_2))x \\
		&=A_1(A_2 x) + B_1(C_2 x) + C_1(A_2 x) + D_1(C_2 x)\\
		&=A_1(\lambda B_2 y) + B_1(\lambda D_2 y)+  C_1(\lambda B_2 y) + D_1(\lambda D_2 y)\\
		&= \lambda ((A_1 B_2 + B_1 D_2)+(C_1 B_2 + D_1 D_2))y.
	\end{align*}
	Hence, $((A_1 A_2 + B_1 C_2)+(C_1 A_2 + D_1 C_2))(\mathcal{B}_M) \subseteq \lambda ((A_1 B_2 + B_1 D_2)+(C_1 B_2 + D_1 D_2))(\mathcal{B}_{M^\perp}).$
	Thus, $T_1T_2 \in \phi_L (M,N,\lambda).$
	
	\vspace{.75cm}  
	
	\noindent (2.) Assume that $T_1 \in \phi_R (M,N,\lambda).$ By similar arguments as in the first case, we can prove 
	$$T_1T_2 \in \phi_R (M,N,\lambda).$$
	
	\noindent (3.) Proof follows from $(1)$ and $(2)$.
	
	\noindent (4.) Proof follows from the fact that $$\phi (M,N,\lambda) \subseteq \psi(M,N,\lambda) \cap \psi (M\p,N\p, \lambda).$$
\end{proof}

For an operator \( T \in \mathcal{B}(\mathcal{H}, \mathcal{K}) \), the analysis of its powers and the series \(\sum\limits_{n=0}^\infty T^n\) is important concept of operator theory. These sequences are pivotal in investigating the spectral characteristics of operators, particularly for compact, self-adjoint, and normal operators. The powers of operators and their series expansions play a significant role in solving operator equations and approximating operators in diverse applications. In quantum mechanics, such series are instrumental in understanding the temporal evolution of quantum states. Furthermore, the study of operator powers is crucial for exploring the asymptotic behavior of operators, the definition of invariant subspaces, and for applications in solving integral equations, as well as enhancing numerical methods for operator approximation.

In this paper, we establish conditions under which the sequence $\{T_n\}$ and series of powers $\sum\limits_{n=0}^{\infty} T^n$ of a complementable operator $T$ remain complementable.

\begin{lemma}\label{lemseries}
	Let $T \in \mathcal{B}(\mathcal{H}).$ Then $\sum\limits_{n=0}^{\infty} \alpha_n T^n$ converges in $\mathcal{B}(\mathcal{H})$  if and only if there exists $\beta <1$ such that $|\alpha_n| ^\frac{1}{n} \|T\| \leq \beta$ for all $n \in \mathbb{N}.$
\end{lemma}
\begin{proof}
	Assume that there exists $\beta < 1$ such that $|\alpha_n|^{\frac{1}{n}} \|T\| \leq \beta$ for all $n \in \mathbb{N}$. 
	This implies 
	$$
	|\alpha_n| \|T\|^n \leq \beta^n \quad \text{for all } n \in \mathbb{N}.
	$$
	Thus, 
	$$
	\sum\limits_{n=0}^{\infty} \|\alpha_n T^n\| \leq \sum\limits_{n=0}^{\infty} |\alpha_n| \|T^n\| \leq \sum\limits_{n=0}^{\infty} \beta^n.
	$$
	Since $\beta < 1$, the geometric series converges, and we have 
	$$
	\sum\limits_{n=0}^{\infty} \|\alpha_n T^n\| \leq \sum\limits_{n=0}^{\infty} \beta^n < \frac{1}{1-\beta}.
	$$
	Since $\mathcal{B}(\mathcal{H})$ is complete and $\sum\limits_{n=0}^{\infty} \alpha_n T^n$ is absolutely summable, it follows that $\sum\limits_{n=0}^{\infty} \alpha_n T^n$ is convergent and $$\|\sum\limits_{n=0}^{\infty} \alpha_n T^n\| \leq \sum\limits_{n=0}^{\infty} \|\alpha_n T^n\| \leq \frac{1}{1-\beta}.$$
	Hence, $$\sum\limits_{n=0}^{\infty} \alpha_n T^n \in \mathcal{B}(\mathcal{H}).$$
	
	Next, suppose there is no $\beta < 1$ such that $|\alpha_n|^{\frac{1}{n}} \|T\| \leq \beta$ for all $n \in \mathbb{N}$. 
	Then, for each $n$, we can find a natural number $n_k \geq n$ such that 
	$$
	1 - \frac{1}{n} \leq |\alpha_{n_k}|^{\frac{1}{n_k}} \|T\|.
	$$
	From this, we deduce 
	$$
	1 \leq \lim\limits_{n_k \to \infty} |\alpha_{n_k}|^{\frac{1}{n_k}} \|T\|.
	$$
	This implies 
	$$
	1 \leq \lim\limits_{n_k \to \infty} |\alpha_{n_k}| \|T\|^{n_k}.
	$$
	Hence, 
	$$
	\lim\limits_{n \to \infty} |\alpha_n| \|T\|^n \neq 0.
	$$
	Thus, the series $\sum\limits_{n=0}^{\infty} \alpha_n T^n$ does not converge.
	
\end{proof}

\begin{theorem}
	Let $T \in \mathcal{B}(\mathcal{H})$ be as in (\ref{T}). Let $T \in \phi(M,N,\lambda).$
	\begin{enumerate}
		\item  $\alpha_nT^n\in \phi(M, N, \lambda)$ for all $n \in \mathbb{N}$.
		
		Moreover, if $\alpha_n T^n \to \overset{\sim}{T} $  uniformly in $\mathcal{B}(\mathcal{H}
		, \mathcal{K})$ and $\overset{\sim}{T} \in \mathcal{C}_D(M,N)$, then $\overset{\sim}{T}\in \psi(M, N, \lambda).$
		\item $S_n=\sum \limits_{i=1}^{n}  \alpha _i T^i\in \psi(M, N, \lambda),$ for all $n \in \mathbb{N}.$
		\item If there exists $\beta <1$ such that $|\alpha_n| ^\frac{1}{n} \|T\| \leq \beta$ for all $n \in \mathbb{N},$ then $\sum\limits_{n=1}^{\infty} \alpha_n T^n$ converges in $\mathcal{B}(\mathcal{H}).$
		Moreover, if $\sum\limits_{n=1}^{\infty} \alpha_n T^n \in \mathcal{C}_D(M,N) $, then $\sum\limits_{n=1}^{\infty} \alpha_n T^n \in \psi(M, N, \lambda).$
	\end{enumerate}
\end{theorem}

\begin{proof}
	We have $\alpha _n T^n = (\alpha _n T^{n-1}) T$. By  Proposition \ref{lem25} (1), we have $$\alpha _n T^n \in \phi_L (M,N,\lambda).$$
	
	Alternatively, we can write $\alpha _n T^n = T (\alpha _n T^{n-1})$. By Proposition \ref{lem25} (2), we have $$\alpha _n T^n \in \phi_R (M,N,\lambda).$$\\
	Combining these results, we conclude that $$\alpha _n T^n \in \phi (M,N,\lambda) \subseteq \psi(M,N,\lambda) \cap \psi (M\p,N\p, \lambda).$$
	
	Now, let us assume that $\alpha _n T^n \to \overset{\sim}{T}$ uniformly. Thus, $\{\alpha _n T^n\}$ is a sequence in $\psi(M,N,\lambda)$ such that $\alpha _n T^n \to \overset{\sim}{T}$ uniformly. By Corollary \ref{boundedlambda}, it follows that $\overset{\sim}{T} \in \psi(M,N,\lambda).$
	
	\noindent  We have $S_n = (\sum \limits_{i=0}^{n} \alpha _{i+1} T^n) T$. By Proposition \ref{lem25} (1), we have $$S_n \in \phi_L(M,N, \lambda).$$
	
	Alternatively, we can write $S_n =  T(\sum \limits_{i=0}^{n} \alpha _{i+1} T^n).$ By Proposition \ref{lem25} (2), we have $$S_n \in \phi_R(M,N, \lambda).$$\\
	Combining these results, we conclude that $$S_n \in \phi (M, N, \lambda)\subseteq \psi(M,N,\lambda) \cap \psi (M\p,N\p, \lambda).$$
	
	Assume that there exists $\beta < 1$ such that $|\alpha_n|^{\frac{1}{n}} \|T\| \leq \beta$ for all $n \in \mathbb{N}$. From \ref{lemseries}, we have 
	$$
	\sum\limits_{n=0}^{\infty} \alpha_n T^n 
	$$ 
	converges in $\mathcal{B}(\mathcal{H})$.\\
	So, 
	$$
	\sum\limits_{n=1}^{\infty} \alpha_n T^n =
	\sum\limits_{n=0}^{\infty} \alpha_n T^n - \alpha_1
	$$ 
	converges in $\mathcal{B}(\mathcal{H})$.\\
	Now, let us take 
	$$
	\sum\limits_{i=1}^{\infty} \alpha_i T^i = T^\#.
	$$
	Assume that $T^\# \in \mathcal{C}_D(M,N)$. Thus, $\{S_n\}$ is a sequence in $\psi(M,N,\lambda)$ such that $S_n \to T^\#$ uniformly. By Corollary \ref{boundedlambda}, it follows that 
	$$
	T^\# \in \psi(M,N,\lambda).
	$$
	
\end{proof}

In contrast to the conditions under which a sequence of complementable operators converge to a complementable operator, it is also possible to approximate operators by sequences of operators that do not necessarily possess this complementability property. This concept is particularly useful for exploring the broader landscape of bounded operators and their approximations, as it provides a framework for understanding the limits of complementability. The following discussion demonstrates how every bounded operator can be realized as the pointwise limit of a sequence of operators that are not necessarily $(M,N)$-complementable from which we can conclude that certain $(M,N)$-complementable are boundary points of non $(M,N)$-complementable operators in the strong operator topology. We denote the set of non $(M,N)$-complementable operators as $(\psi(M,N))^\complement$.

\begin{theorem}\label{ptwise}
	Let $ T \in \psi(M,N)$ be as defined in \ref{T}. If $ \mathcal{R}(C) $ is infinite-dimensional, then $ T \in \partial  ((\psi(M,N))^\complement).$
\end{theorem}

\begin{proof}
	Consider the operator 
	$$
	T = \begin{pmatrix}
		A & B \\
		C & D
	\end{pmatrix}.
	$$ 
	Given that $T$ is $(M,N)$-complementable, it follows that $\mathcal{R}(C) \subseteq \mathcal{R}(D)$.\\
	Now, $\mathcal{R}(C)$ is infinite dimensional, this gives $\mathcal{R}(D)$ is also infinite dimensional. Now consider $\overset{\sim}{D} = D|_{\mathcal{N}(D)^\perp} : \mathcal{N}(D)^\perp \to \mathcal{R}(D)$ defined by $\overset{\sim}{D}(x) = D(x).$ Thus $\overset{\sim}{D}$ is bijective. 
	
	Now, observe that 
	$$ 
	D^{-1}(\mathcal{R}(C)) \cap \mathcal{N}(D)^\perp = (\overset{\sim}{D})^{-1}(\mathcal{R}(C))
	$$ 
	is infinite-dimensional as $\mathcal{R}(C)$ is infinite-dimensional and $\overset{\sim}{D}$ is bijective. 
	
	Let $\mathcal{A} = \{a_1, a_2, \ldots\}$ be a countably infinite orthonormal set in $D^{-1}(\mathcal{R}(C)) \cap \mathcal{N}(D)^\perp$. Notably, $0 \neq D(a_i) \in \mathcal{R}(C)$ for all $i$, and $D(a_i) \neq D(a_j)$ for all $i \neq j$. 
	
	Let $\mathcal{A}_n = \{a_n, a_{n+1}, \ldots\}$.
	\\
	Define $D_n: M^\perp \to N^\perp$ by
	$$
	D_n(x) = 
	\begin{cases} 
		D(x), & \text{if } x \in (\text{span}({\mathcal{A}}_n))^\perp \cap \mathcal{N}(D)^\perp, \\
		0, & \text{if } x \in (\overline{\text{span}({\mathcal{A}}_n)} \cap \mathcal{N}(D)^\perp) \oplus \mathcal{N}(D).
	\end{cases}
	$$
	Now, consider a sequence of operators 
	$$
	T_n = \begin{pmatrix}
		A & B \\
		C & D_n
	\end{pmatrix}.
	$$ 
	Clearly, $0 \neq b_n = D(a_n) \in \mathcal{R}(C)$ for all $n$. Suppose $b_n \in \mathcal{R}(D_n)$; then $b_n = D(c_n)$ for some $c_n \in (\text{span}({\mathcal{A}}_n))^\perp \cap \mathcal{N}(D)^\perp$.
	Therefore, $$a_n - c_n \in \mathcal{N}(D) \cap \mathcal{N}(D)^\perp.$$ Thus, $$a_n = c_n.$$\\
	However, by construction, $D_n(a_n) = 0$, so $D_n(c_n) = 0$, which leads to a contradiction.\\
	Consequently, $b_n \notin \mathcal{R}(D_n),$ implying $\mathcal{R}(C) \nsubseteq \mathcal{R}(D_n).$
	
	Next, we prove that $\{T_n\}$ converges to strongly $T$. To prove this, it is sufficient to prove $\{D_n\}$ converges strongly to $D$. \\
	\textbf{Case 1:} If $x \in \mathcal{N}(D)$, then $D(x) = 0 = D_n(x)$ for all $n$.\\
	\textbf{Case 2:} If $x \in (\text{span}({\mathcal{A}}))^\perp \cap \mathcal{N}(D)^\perp$, then $x \in (\text{span}({\mathcal{A}}_n))^\perp \cap \mathcal{N}(D)^\perp$ for all $n$, so $$D_n(x) = D(x) \text{ for all } n.$$
	
	\noindent \textbf{Case 3:} Let $x \in \overline{\text{span}({\mathcal{A}})} \cap \mathcal{N}(D)^\perp.$ Let $\{x_n\}$ be a sequence in  $ \text{span}({\mathcal{A}}) \cap \mathcal{N}(D)^\perp$ such that $x_n \to x.$
	
	Let $x = \sum\limits_{i=1}^{\infty} \alpha_i a_i$, with $x_i \in {\mathcal{A}}$. Now $D_n(x)= \sum\limits_{i=1}^{n-1} \alpha_i D_n(a_i)+ \sum\limits_{i=n}^{\infty} \alpha_i D_n(a_i)=\sum\limits_{i=1}^{n-1} \alpha_i D(a_i).$ So we have $$\lim_{n \to \infty} D_n x= \sum\limits_{i=1}^{\infty} \alpha_i D(a_i).$$
	Also, \begin{align*}
		Dx =& D(\sum\limits_{i=1}^{\infty} \alpha_i a_i) \\
		&= D(\lim_{n \to \infty} \sum\limits_{i=1}^{n} \alpha_i a_i)\\
		&=\lim_{n \to \infty}D( \sum\limits_{i=1}^{n} \alpha_i a_i)\\
		&= \lim_{n \to \infty}( \sum\limits_{i=1}^{n} \alpha_i D(a_i))\\
		&=\sum\limits_{i=1}^{\infty} \alpha_i D(a_i).
	\end{align*} 
	So, we have $$\lim_{n\to \infty}D_n x = Dx.$$
	Thus, $$ T \in \partial ( (\psi(M,N))^\complement).$$
	
\end{proof}
\begin{theorem}\label{ptwise2}
	Let $T \in \psi(M,N)$ be as defined in \ref{T}. If $\mathcal{R}(C)$ is finite-dimensional and $\mathcal{R}(D)$ is infinite-dimensional, then $T \in \partial ( (\psi(M,N))^\complement )$.
\end{theorem}

\begin{proof}
	Consider the operator 
	$$
	T = \begin{pmatrix}
		A & B \\
		C & D
	\end{pmatrix}.
	$$
	Since $T$ is $(M, N)$-complementable, we have $\mathcal{R}(C) \subseteq \mathcal{R}(D)$.
	
	Now consider $\overset{\sim}{D} = D|_{\mathcal{N}(D)^\perp} : \mathcal{N}(D)^\perp \to \mathcal{R}(D)$ defined by $\overset{\sim}{D}(x) = D(x).$ Now $\overset{\sim}{D}$ is bijective. 
	
	Observe that 
	$$ 
	D^{-1}(\mathcal{R}(C)) \cap \mathcal{N}(D)^\perp = (\overset{\sim}{D})^{-1}(\mathcal{R}(C))
	$$ 
	is finite-dimensional as $\mathcal{R}(C)$ is finite-dimensional and $\overset{\sim}{D}$ is bijective, while 
	$$ 
	\mathcal{N}(D)^\perp = (\overset{\sim}{D})^{-1}(\mathcal{R}(D))
	$$ 
	is infinite-dimensional as $\mathcal{R}(D)$ is infinite-dimensional and $\overset{\sim}{D}$ is bijective.
	
	Let $\mathcal{A} = \{a_1, a_2, a_3, a_4, \ldots\}$ be a countably infinite orthonormal set in  
	$$ 
	\mathcal{N}(D)^\perp \setminus \big(D^{-1}(\mathcal{R}(C)) \cap \mathcal{N}(D)^\perp\big). 
	$$  
	Note that $0 \neq D(a_i) \notin \mathcal{R}(C)$ for all $i$, and $D(a_i) \neq D(a_j)$ for $i \neq j$.  
	Let $\mathcal{A}_n = \{a_n, a_{n+1}, \ldots\}$.\\
	Define $D_n: M^\perp \to N^\perp$ by
	$$
	D_n(x) = 
	\begin{cases} 
		D(x), & \text{if } x \in \text{span}({\mathcal{A}}_n)^\perp \cap \mathcal{N}(D)^\perp, \\
		0, & \text{if } x \in (\overline{\text{span}({\mathcal{A}}_n)} \cap \mathcal{N}(D)^\perp)\oplus \mathcal{N}(D).
	\end{cases}
	$$
	Since $\mathcal{R}(C)$ is finite-dimensional, $\mathcal{N}(C)^\perp$ is finite-dimensional. Let $\mathcal{F} = \{f_1, f_2, \ldots\}$ be a countable orthonormal set in $\mathcal{N}(C)$.\\
	Now, Define $h: \text{span}(\mathcal{F}) \to \text{span}({\mathcal{A}})$ as follows  
	$$h (f_n) = a_n, \text{ for all } f_n.$$ From this, we can extend $h$ to $\text{span}(\mathcal{F}).$ Note that $h : \text{span}(\mathcal{F}) \to \text{span}({\mathcal{A}})$ is a bounded operator. So $h$ has a bounded extension to $\overline{\text{span}(\mathcal{F})}.$ 
	
	\noindent Define $C_n: M \to N^\perp$ by  
	$$
	C_n(x) = 
	\begin{cases} 
		C(x), & \text{if } x \in \mathcal{N}(C)^\perp, \\
		D(h(x)), & \text{if } x \in \overline{\text{span}(\mathcal{F}_n)}, \\
		0, & \text{if } x \in (\text{span}(\mathcal{F}_n))^\perp \cap \mathcal{N}(C).
	\end{cases}
	$$
	Now, consider the sequence of operators  
	$$
	T_n = \begin{pmatrix}
		A & B \\
		C_n & D_n
	\end{pmatrix}.
	$$
	
	For all $n$, we have $f_n \in \text{span}(\mathcal{F}_n)$, so $C_n(f_n) = D(h(f_n)) = D(a_n).$  
	Let $0 \neq b_n = D(a_n) \in \mathcal{R}(C_n)$. Suppose $b_n \in \mathcal{R}(D_n)$. Then $b_n = D(c_n)$ for some $c_n \in (\text{span}({\mathcal{A}}_n))^\perp \cap \mathcal{N}(D)^\perp$. Hence,  
	$$
	a_n - c_n \in \mathcal{N}(D) \cap \mathcal{N}(D)^\perp,
	$$
	which implies $a_n = c_n$. However, by construction, $D_n(a_n) = 0$, leading to $D_n(c_n) = 0$, a contradiction. Thus, $b_n \notin \mathcal{R}(D_n)$. Since $b_n = D(a_n) = C_n(f_n)$, we conclude that $\mathcal{R}(C_n) \nsubseteq \mathcal{R}(D_n)$.
	
	To complete the proof, we need to show that $T_n \to T$ strongly. It suffices to prove $D_n \to D$ strongly and $C_n \to C$ strongly.\\
	
	\noindent \textbf{Strong convergence of $\{D_n\}$:}  
	From the proof of Theorem \ref{ptwise}, we have $\{D_n\}$ converges strongly to $D$.\\
	
	\noindent \textbf{Strong convergence of $\{C_n\}$:}  \\
	\noindent \textbf{Case 1:} If $x \in \mathcal{N}(C)^\perp$, then $C_n(x) = C(x)$ for all $n$. Hence,  
	$$C_n x \to Cx.$$

	\noindent \textbf{Case 2:} If $x \in \overline{\text{span}(\mathcal{F})} \cap \mathcal{N}(C)$, then $x = \sum_{i=1}^{\infty} \alpha_i f_i$, with $f_i \in \mathcal{F}$. Let  
	$$x_n = \sum_{i=1}^{n} \alpha_i f_i.$$  
	Then $x_n \to x$. 
	
	Now $C_n(x)= \sum\limits_{i=1}^{n-1} \alpha_i C_n(f_i)+ \sum\limits_{i=n}^{\infty} \alpha_i C_n(f_i)=\sum\limits_{i=n}^{\infty} \alpha_i D(h(f_i))=\sum\limits_{i=n}^{\infty} \alpha_i D(a_i).$ So we have $$\lim_{n \to \infty} C_n x= \lim\limits_{n\to \infty} \sum\limits_{i=n}^{\infty} \alpha_i D(a_i).$$
	
	Since $h$ is bounded and $x_n=\sum \limits_{i=1}^{n} \alpha_i f_i \to x$, we have $h(x_n) \to h(x).$ So, $\sum \limits_{i=1}^{n} \alpha_i a_i \to h(x)=\sum \limits_{i=1}^{\infty} \alpha_i a_i.$ Since $D$ is bounded, we have $D\Big(\sum \limits_{i=1}^{n} \alpha_i a_i \Big) \to D\Big(\sum \limits_{i=1}^{\infty} \alpha_i a_i\Big).$ This gives $$\Big \|\sum \limits_{i=n}^{\infty} \alpha_i D(a_i)\Big\| = \Big\|\sum \limits_{i=1}^{\infty} \alpha_i D(a_i)- \sum \limits_{i=1}^{n-1} \alpha_i D(a_i)\Big\| \to 0 .$$
	Thus we have $\lim_{n \to \infty} C_n x = 0.$ \\
	Also, $x \in \mathcal{N}(C),$ we have $Cx=0.$\\
	Thus $$\lim_{n \to \infty} C_n x = Cx.$$

	\noindent \textbf{Case 3:} If $x \in (\text{span}(\mathcal{F}))^\perp \cap \mathcal{N}(C)$, then $x \in (\text{span}(\mathcal{F}))^\perp \subseteq (\text{span}(\mathcal{F}_n))^\perp$ for all $n$, so $C_n(x) = 0$ for all $n$. Also, $x \in \mathcal{N}(C)$, so $Cx = 0$, and $C_n x \to Cx$.  
	
	\noindent Thus, $C_n \to C$ strongly, and $T_n \to T$ strongly. Therefore, $$T \in \partial((\psi(M,N))^\complement).$$
\end{proof}

\begin{corollary}
	Let $T \in \psi(M,N).$ If $T \notin \partial((\psi(M,N))^\complement)$, then $\mathcal{R}(D)$ is finite-dimensional.
\end{corollary}
\begin{proof}
	From Theorem \ref{ptwise} and Theorem \ref{ptwise2}, we have $T \in \partial((\psi(M,N))^\complement),$ if $\mathcal{R}(D).$ From this we conclude that If $T \notin \partial((\psi(M,N))^\complement)$, then $\mathcal{R}(D)$ is finite-dimensional.
\end{proof}

\section{Topological Properties of the Set of Complementable Operators}
In this section, we derive results concerning the topological properties of certain subsets of complementable operators. These operators satisfying specific structural conditions form a significant class in operator theory. By exploring their topological characteristics, we aim to uncover insights into their behaviour with respect to convergence, continuity, and other related properties within the context of operator norms. 

\begin{theorem}\label{Closure}
	$\overline{\psi(M,N,\lambda)} \cap \mathcal{B}\mathcal{C}_D(M,N) \subseteq \psi(M,N,\lambda).$
\end{theorem}
\begin{proof}
	Let $T = \begin{pmatrix} A & B \\ C & D \end{pmatrix} \in \overline{\psi(M,N,\lambda)} \cap \mathcal{B}\mathcal{C}_D(M,N)$. Let $T_n = \begin{pmatrix} A_n & B_n \\ C_n & D_n \end{pmatrix}$ be a sequence in $\psi(M,N,\lambda)$ such that $T_n \to T$. Since $\mathcal{R}(D)$ is closed, by applying Theorem \ref{conv}, we conclude that $T \in \psi(M,N)$.
	
	Now we can write $Cx = \lim C_n x = \lim \lambda D_n y_n$, with $\|y_n\| \leq 1$. Since $T \in \psi(M,N)$, we have $Cx = Dz$ for some $z \in \mathcal{B}_{(\mathcal{N}(D))\p}$. Furthermore, we have 
	$$
	\|\lambda D y_n - \beta_n D_n y_n\| \leq \lambda \|D - D_n\| \|y_n\| \leq \lambda \|D - D_n\| \to 0.
	$$
	This implies that $\{\lambda D y_n\}$ converges to the same limit to which $\{\lambda D_n y_n\}$ converges. Thus, $\lambda D y_n \to Dz$.
	So we have 
	$$
	z = D\+ Dz = \lambda \lim_{n\to \infty} D\+ D y_n.
	$$
	Therefore, $\|z\| = \lim\limits_{n\to \infty} \|\lambda D\+ Dy_n\| \leq \lim \limits_{n\to \infty}\|\lambda y_n\| \leq \lambda$.
	Consequently, 
	$$
	Cx = \lambda D\left(\frac{z}{\lambda}\right).
	$$
	Hence, $T \in \psi(M,N,\lambda)$.
\end{proof}

\begin{corollary}
	$\psi(M,N,\lambda)$ is closed in $\mathcal{B}\mathcal{C}_D(M,N).$
	
\end{corollary}
\begin{proof}
	From Theorem \ref{Closure}, we have 
	$$
	\overline{\psi(M,N,\lambda)} \cap \mathcal{B}\mathcal{C}_D(M,N) \subseteq \psi(M,N,\lambda).
	$$
	Hence 
	$$
	\overline{\psi(M,N,\lambda)} \cap \mathcal{B}\mathcal{C}_D(M,N) \subseteq \psi(M,N,\lambda) \cap \mathcal{B}\mathcal{C}_D(M,N).
	$$
	On the other hand, it is easy to see that $\psi(M,N,\lambda) \subseteq \overline{\psi(M,N,\lambda)}$, which leads to 
	$$
	\psi(M,N,\lambda) \cap \mathcal{B}\mathcal{C}_D(M,N) \subseteq \overline{\psi(M,N,\lambda)} \cap \mathcal{B}\mathcal{C}_D(M,N).
	$$
	Hence, $$\overline{\psi(M,N,\lambda)} \cap \mathcal{B}\mathcal{C}_D(M,N) = \psi(M,N,\lambda) \cap \mathcal{B}\mathcal{C}_D(M,N).$$
	Hence, $\psi(M,N,\lambda)$ is closed in $\mathcal{B}\mathcal{C}_D(M,N).$
\end{proof}

\begin{corollary}
	Let $\mathcal{K}$ be a finite-dimensional Hilbert space. Then $\psi(M,N,\lambda)$ is a closed subset of $\mathcal{B}(\mathcal{H}, \mathcal{K})$.
\end{corollary}
\begin{proof}
	We have $\mathcal{B}\mathcal{C}_D(M,N)= \{ T \in \mathcal{B}(\mathcal{H},\mathcal{K}) : \mathcal{R}(D) \text{ is closed}\}.$ When $\mathcal{K}$ is a finite-dimensional Hilbert space,  $\mathcal{R}(D) $ is always closed, for any $T \in \mathcal{B}(\mathcal{H},\mathcal{K}).$
	Therefore we have $$\mathcal{B}\mathcal{C}_D(M,N)=\mathcal{B}(\mathcal{H},\mathcal{K}).$$
	By Theorem \ref{Closure}, we have $$\overline{\psi(M,N,\lambda)} \cap \mathcal{B}\mathcal{C}_D(M,N) \subseteq \psi(M,N,\lambda).$$
	Hence $$\overline{\psi(M,N,\lambda)}  \subseteq \psi(M,N,\lambda).$$
\end{proof}

\begin{center}
	{\bf Acknowledgements }
\end{center}
The first author thanks the National Institute of Technology Karnataka (NITK), Surathkal for the financial support. The present work of the second author was partially supported by Science and Engineering Research Board (SERB), Department of Science and Technology, Government of India (Reference Number: MTR/2023/000471) under the scheme ``Mathematical Research Impact Centric Support (MATRICS)".

\bibliographystyle{plain}

\end{document}